\documentclass[oneside, 11pt]{amsart} 
\usepackage[margin=1in]{geometry}
\usepackage{amsmath, amssymb, amsfonts, amstext, amsthm, amscd}
\usepackage[mathscr]{euscript}
\usepackage{enumerate}
\usepackage{tikz,color,soul}
\usepackage[english]{babel}
\usepackage{chngcntr}
\usetikzlibrary{arrows.meta, positioning,arrows}
\newtheorem{theorem}{Theorem}[section]

\newtheorem{prop}[theorem]{Proposition}
\newtheorem{lem}[theorem]{Lemma}
\newtheorem{cor}[theorem]{Corollary}
\theoremstyle{definition}
\newtheorem{definition}[theorem]{Definition}
\newtheorem{eg}[theorem]{Example}
\newtheorem{remark}[theorem]{Remark}

\numberwithin{equation}{section}
\usepackage{tikz}
\usetikzlibrary{arrows.meta}

\title{On a McCoy-like condition for rings}
\author{Sharang Thimmaiah}
\author{Raisa DSouza}
\date{}
\subjclass{16U80; 16S15}
\begin{document}

\maketitle
\begin{abstract}
    We study rings $R$ for which whenever non-zero polynomials $f(x)$ and $g(x)$  satisfy $f(x)g(x)f(x)=0$, it implies that there is a non-zero element $r\in R$ such that $f(x)rf(x)=0$. We call such rings inner McCoy rings. We explore some examples of rings that are inner McCoy, determine relationships between the class of inner McCoy rings and some known classes of rings. Furthermore, we show that the maximal inner McCoy subring of a matrix ring is the triangular matrix ring. Finally, we construct new inner McCoy rings from known ones. 
    
    \vspace{0.5cm}
    \noindent
    \textbf{Keywords:} Inner annihilators; McCoy rings; Inner McCoy rings; Maximal subring.
\end{abstract}

\section{Introduction}
All rings considered in this paper will be associative with unity.

In 1942, Neal H. McCoy proved a result that has become commonplace in most introductory commutative algebra textbooks. This relates the nonzero zero divisors of the polynomial ring of a commutative ring with unity, with a nonzero annihilator taken from the base ring:
\begin{lem}\cite[Theorem 2]{mccoy1942}
    Suppose $R$ is a commutative ring with unity. Let $f(x)$ and $g(x)$ be two nonzero elements of the ring $R[x]$ such that $f(x)g(x)=0$. Then, we can find nonzero elements $s, t \in R$ such that $sg(x) = 0$ and $f(x)t = 0$.
\end{lem}

This theorem has been rephrased in literature by looking at $R[x]$ as a $R$-module and attempting to find a nonzero right (left) annihilator of $f(x)$ (or $g(x)$) from $R$, provided a nonzero right (left) annihilator exists for $f(x)$ (or $g(x)$) in $R[x]$. For examples one can see \cite{hong2000} and \cite{hong2010}.

The natural question to ask after this was if a similar result held for rings that were not necessarily commutative. This was answered in the negative by Jones and Weiner \cite{weiner1952} with a noncommutative example. Rings that do follow such a property, without necessarily being commutative, were first studied by Rege and Chhawchharia \cite{rege1997armendariz} and were called ``McCoy rings". We define them below:
\begin{definition}\cite{rege1997armendariz,nielsen2006}
    Let $R$ be an associative ring with unity. We call $R$ a \emph{left McCoy ring} if for every non-zero pair $f(x),g(x) \in R[x]$ such that $f(x)g(x) = 0$ we have a nonzero $s \in R$ such that $sg(x) = 0$. \emph{Right McCoy rings} are defined analogously and rings that are both left McCoy and right McCoy are termed \emph{McCoy rings}.
\end{definition}

Nielsen \cite{nielsen2006} studied the connections between McCoy rings and semi-commutative rings. 

In 2014, Lam and Nielsen \cite{lam2014} introduced a variant of the traditional annihilator, which they termed an ``inner annihilator.'' They introduced this in the context of studying generalisations of Jacobson's lemma using Drazin inverses. We reproduce their definition below:
\begin{definition}\cite{lam2014}
    Let $R$ be an associative ring and let $\alpha \in R$. The \emph{inner annihilators} of $\alpha$, denoted as $\operatorname{IAnn}_R(\alpha)$, is the set of all elements $p \in R$ that annihilate $\alpha$ from the inside, that is, $\alpha p\alpha = 0$.
\end{definition}
Lam and Nielsen also explored some rather basic properties of this new structure, like the fact that the left and right annihilators of $\alpha$ are contained in $\operatorname{IAnn}_R(\alpha)$ and how the set of inner inverses of $\alpha$ is related to $\operatorname{IAnn}_R(\alpha)$.

Apart from these, some other well known classes of rings include \\

   \begin{tabular}{rl}
    reversible: & A ring $R$ in which $ab=0\implies ba=0$, for all $a,b\in R$.\\
    semi-commutative:& A ring $R$ in which $ab = 0\implies aRb = 0$, for all $a,b\in R$.\\
    \emph{i}-reversible:& A ring $R$ in which $ab$ non-zero idempotent $\implies ba$ idempotent, for all $a,b\in R$. \\
\end{tabular}
\\

Rings that are \textit{i}-reversible are also known in literature as ``quasi-reversible'' rings, following work by Jung et al \cite{jung2019}.
\\
Some known results in literature include:
\begin{enumerate}
    \item Reversible rings are McCoy, semi-commutative and \textit{i}-reversible. %This follows from the definition \cite{khurana2020}.
   % \item Reversible rings are semi-commutative.
    \item There is a semi-commutative ring that is not McCoy and there is a McCoy ring that is not semi-commutative \cite{nielsen2006, zhao2009}.
    %\item McCoy rings are right/left McCoy. This follows from the definition.
    \item There are right/left McCoy rings that are not McCoy \cite{nielsen2006}.
\end{enumerate}

%The known relationships between these classes of rings are given in the diagram below:

%\begin{figure}[h]
%    \centering
%\tikzstyle{ringlabel} = [rectangle, 
%minimum width=3cm, 
%minimum height=1cm,
%text centered, 
%draw=white, 
%fill=white]
%\begin{tikzpicture}[node distance=1cm]

%\node[draw] at (0,-3) (rev)[ringlabel] {Reversible};

%\node[draw] at (-5,0) (rmc)[ringlabel] {Right McCoy};

%\node[draw] at (-5,-3) (mc)[ringlabel] {McCoy};

%\node[draw] at (5,-3) (irev)[ringlabel] {\textit{i}-Reversible};

%\draw [-Implies, line width=1pt,double distance=2pt] (rev) -- (semicom);

%\draw [-Implies, line width=1pt,double distance=2pt] (rev) -- (rmc);

%\draw [-Implies, line width=1pt,double distance=2pt] (rev) -- (mc);

%\draw (semicom) -- (rmc) node[midway, fill=white]{\cite{nielsen2006}};

%\draw [-Implies, line width=1pt,double distance=2pt] (rev) -- (irev);

%\draw [-Implies, line width=1pt,double distance=2pt] (mc) -- (rmc);

%\draw (semicom) -- (mc) node [near start, fill=white] {\cite{nielsen2006, zhao2009}};
%\end{tikzpicture}
%\caption{Relationships between known classes of rings}
%    \label{fig1}
%\end{figure}

In this article, we introduce a new ring, which we will call an ``inner McCoy ring''. Apart from studying the relationship between this new class of rings and the existing classes of rings (mentioned above), we also explore some ring extensions of inner McCoy rings and construct new inner McCoy rings from old ones.

\section{Inner McCoy Rings}
Inspired by the two seemingly unrelated concepts of McCoy rings and the inner annihilator, we define a ring to be inner McCoy as follows:
\begin{definition}
    Let $R$ be a ring. We say that $R$ is \emph{inner McCoy} if for every pair $f(x), g(x) \in R[x]\setminus \{0\}$ such that $f(x)g(x)f(x) = 0$, then there exists $r \in R\setminus\{0\}$ such that $f(x)rf(x) = 0$. 
\end{definition}
This definition was motivated by replacing the right (left) annihilator condition in the definition of right (left) McCoy rings as defined by Nielsen, with the inner annihilator as defined by Lam and Nielsen.

Now, as McCoy and inner McCoy rings are similarly defined, we would like to establish a relationship between both classes of rings. This we establish below:

\begin{theorem}\label{MisIM}
    Every McCoy ring is inner McCoy.
\end{theorem}
\begin{proof}
    Suppose that $R$ is McCoy. Suppose for non-zero polynomials $f(x)$ and $g(x)$ in $R[x]$, we have 
$f(x)g(x)f(x) = 0$.
    If $g(x)f(x) = 0$, as $R$ is McCoy, it is left McCoy, so we have a nonzero $r \in R$ such that, 
    $
        rf(x) = 0\Rightarrow f(x)\cdot r \cdot f(x) = 0
    $
    which proves the inner McCoy condition.

    If instead $g(x)f(x) \neq 0$. Consider $f(x)(g(x)f(x)) = 0$. As $R$ is McCoy and hence right McCoy, we have a non-zero $r \in R$, such that,
    $
        f(x)r = 0\Rightarrow f(x)\cdot r \cdot f(x) = 0
    $  
    which again proves the inner McCoy condition.
Therefore in all cases, $R$ must be inner McCoy also.
\end{proof}

 Nielsen \cite{nielsen2006} showed that all reversible rings are McCoy. In view of this and Theorem \ref{MisIM}, all reversible rings are  inner McCoy.

Having shown that all McCoy rings are inner McCoy, we wish to find a non-trivial example of an inner McCoy ring. We provide one such example below.

\begin{theorem}\label{x1}
    Let $R$ be a commutative ring with unity.  Then the ring of all $2\times2$ upper triangular matrices with entries from $R$, namely $T_2(R)$, is inner McCoy.
\end{theorem}
\begin{proof}
    Suppose we have $F(x), G(x) \in T_2(R)[x]\setminus\{0\}$ such that $F(x)G(x)F(x) = 0$. Suppose also that these elements are of the form
    $$ F(x) = \begin{bmatrix}
        f_{11}(x) & f_{12}(x) \\
        0 & f_{22}(x)
    \end{bmatrix}
    , G(x) = \begin{bmatrix}
        g_{11}(x) & g_{12}(x) \\
        0 & g_{22}(x)
    \end{bmatrix}
    $$
    Then the above expression may be rewritten as 
    $$
    \begin{bmatrix}
        f_{11}(x)g_{11}(x)f_{11}(x) & f_{11}(x)g_{11}(x)f_{12}(x) + f_{11}(x)g_{12}(x)f_{22}(x) + f_{12}(x)g_{22}(x)f_{22}(x) \\
        0 & f_{22}(x)g_{22}(x)f_{22}(x)
    \end{bmatrix} = \begin{bmatrix}
        0 & 0 \\
        0 & 0
    \end{bmatrix}.
    $$
    It is enough to find a nonzero matrix $A = \begin{bmatrix}
        a & b \\
        0 & c
    \end{bmatrix}$ such that $F(x)AF(x) = 0$.
    We have the following cases:
    \begin{enumerate}
        \item Suppose $f_{11}(x) = 0$. Then clearly $F(x)E_{11}F(x) = 0$. Similarly if $f_{22}(x) = 0$, then $F(x)E_{22}F(x) = 0$.
        \item Suppose $f_{11}(x) \neq 0 \neq f_{22}(x)$ and without loss of generality $g_{11}(x) \neq 0$. From the expression for $F(x)G(x)F(x) = 0$, we get $f_{11}(x)g_{11}(x)f_{11}(x) = 0$. Suppose $f_{11}(x)g_{11}(x) = 0$. Then as $R$ is commutative and hence McCoy, we have a nonzero $r \in R$ such that $f_{11}(x)r = 0$. Then $F(x)(rE_{11})F(x) = 0$. If instead $f_{11}g_{11} \neq 0$, then as $(f_{11}(x)g_{11}(x))f_{11}(x) = 0$, we have a nonzero $r\in R$ such that $rf_{11}(x) = 0$. As $R$ is commutative we have $f_{11}r=0$. Here again $F(x)(rE_{12})F(x) = 0$.
        \item Suppose $f_{11}(x) \neq 0 \neq f_{22}(x)$ and $g_{11}(x) = 0 = g_{22}(x)$. Then to ensure that $G(x)$ is nonzero, we need $g_{12}(x) \neq 0$. From the expression for $F(x)G(x)F(x) = 0$, we get $f_{11}(x)g_{12}(x)f_{22}(x) = 0$. Suppose $f_{11}(x)g_{12}(x) = 0$. Then as $R$ is also McCoy, we have a nonzero $r \in R$ such that $f_{11}(x)r = 0$. Then $F(x)(rE_{12})F(x) = 0$. If instead $f_{11}g_{12} \neq 0$, then as $(f_{11}(x)g_{12}(x))f_{22}(x) = 0$, we have a nonzero $r\in R$ such that $rf_{22}(x) = 0$. Here again $F(x)(rE_{12})F(x) = 0$.
    \end{enumerate}
    This exhausts all cases and hence $T_2(R)$ is inner McCoy.
\end{proof}

\begin{remark}
It is well known in literature that upper triangular matrix rings over any nonzero ring of any size $n \geq 2$ is never McCoy \cite[Prop. 10.2]{camillo2008}. Example \ref{x1} shows that there are examples of rings that are not McCoy but are inner McCoy, further showing that the class of inner McCoy rings is strictly larger than that of McCoy rings. 
\end{remark}

\begin{remark}
    One can show by induction that the triangular matrix ring $T_n(R)$ is inner McCoy when $R$ is a commutative ring.
\end{remark}

Now, we present an example of a ring that is not inner McCoy.%, to ensure that this is not a class of rings that contains all possible rings (that is, it is not a property common to all rings, and so is not particularly interesting).

\begin{eg}\label{mn}
    Consider the matrix ring $M_2(\mathbb{R})$, consisting of all $2\times 2$ matrices with real valued entries. $M_2(\mathbb{R})$ is not inner McCoy.

    We first note that any polynomial from $M_2(\mathbb{R})[x]$ may be rewritten as a $2\times 2$ matrix whose entries are taken from the polynomial ring $\mathbb{R}[x]$. Consider the polynomials
    $$A(x) = \begin{bmatrix}
        x & x^2 \\
        x^3 & x^4
    \end{bmatrix}, B(x) = \begin{bmatrix}
        x & -x^2 \\
        -1 & x
    \end{bmatrix} \in M_2(\mathbb{R})[x]\setminus\{0\}.$$
    Clearly $A(x)B(x) = 0$, which means $A(x)B(x)A(x) = 0$. 
    If $D = \begin{bmatrix}
        a & b \\
        c & d
    \end{bmatrix}$ is such that $A(x)DA(x) = 0$, then $$
    \begin{bmatrix}
        x & x^2 \\
        x^3 & x^4
    \end{bmatrix}
    \begin{bmatrix}
        a & b \\
        c & d
    \end{bmatrix}
    \begin{bmatrix}
        x & x^2 \\
        x^3 & x^4
    \end{bmatrix} = \begin{bmatrix}
        0 & 0 \\
        0 & 0
    \end{bmatrix}.
    $$ Solving this gives us $a = b = c = d = 0$, that is, $D=0$. Thus $M_2(\mathbb{R})$ is not inner McCoy.
\end{eg}

Here, we provide for a large class of rings that are not inner McCoy.

\begin{eg}
    Let $R$ be a commutative ring with unity. For all $n \geq 2$, the matrix ring $M_n(R)$ is not inner McCoy. To show this, just consider the polynomials of $M_n(R)[x]:$ $$ A(x) = \begin{bmatrix}
        x & x^2 & \cdots & x^n \\
        x^{n+1} & x^{n+2} & \cdots & x^{2n} \\
        \vdots & \vdots & \ddots & \vdots\\
        x^{n^2-n+1} & x^{n^2-n+2} & \cdots & x^{n^2}
    \end{bmatrix}, B(x) = \begin{bmatrix}
        -x^{n-1} & -x^{n-1} & \cdots & -x^{n-1} \\
        x^{n-2} & x^{n-2} & \cdots & x^{n-2} \\
        \vdots & \vdots & \ddots & \vdots \\
        1 & 1 & \cdots & 1
    \end{bmatrix}.$$
    It is clear that $A(x)B(x)A(x) = 0$, but solving for a $D \in M_n(R)$ such that $A(x)DA(x) = 0$ gives us $D=0$ always.
\end{eg}

\section{Semi-commutative, \emph{i}-reversible and Inner McCoy Rings}\label{rel}
In this section we first, discuss the relationship between inner McCoy rings and  \emph{i}-reversible rings. These discussions involve the diagonal upper triangular matrix ring, $DUT_n(R)$.
 Here $DUT_n(R)$ is the ring of upper triangular matrices of order $n$ whose entries come from the ring $R$ and with all diagonal entries equal. 

We recall some useful results about the diagonal upper triangular matrix ring proved by Lei, Chen and Ying \cite{lei2007} 
and  Lama et al. \cite{lama2022}.
\begin{prop}\cite[Theorem 2]{lei2007} \label{p1}
 $R$ is a McCoy ring if and only if the ring $DUT_n(R)$ is McCoy for all $n \geq 2$. 
 \end{prop}

\begin{prop}\cite[Theorem 3.1]{lama2022} \label{c2}
    Let $R$ be a ring and let $n \geq 2$ be a positive integer. The ring $DUT_n(R)$ is \emph{i}-reversible if and only if $R$ has only trivial idempotents.
\end{prop}

The reason we bring these results is to make a few observations based on this ring $DUT_n(R)$ and its relation between the inner McCoy property and \emph{i}-reversibility.

An analogue of Proposition \ref{p1} holds for inner McCoy rings. We state and prove it below:

\begin{theorem}\label{c1}
    $R$ is an inner McCoy ring if and only if the ring $DUT_n(R)$ is inner McCoy where $n \geq 2$.
\end{theorem}
\begin{proof}
Let $R$ be inner McCoy and let $G(x), F(x) \in DUT_n(R)[x]\setminus\{0\}$ such that $F(x)G(x)F(x) = 0$. Here, we suppose that 
\[G(x) = 
\begin{bmatrix}
    g(x)    &   g_{12}(x)   &   g_{13}(x)   &   \cdots  &   g_{1n}(x)\\
    0   &   g(x)    &   g_{23}(x)   &   \cdots  &   g_{2n}(x)\\
    0   &   0   &   g(x)    &   \cdots  &   g_{3n}(x)\\
    \vdots  &   \vdots  &   \vdots  &   \ddots  &   \vdots\\
    0   &   0   &   0   &   \cdots  &   g(x)\\
\end{bmatrix},\quad
F(x) = 
\begin{bmatrix}
    f(x)    &   f_{12}(x)   &   f_{13}(x)   &   \cdots  &   f_{1n}(x)\\
    0   &   f(x)    &   f_{23}(x)   &   \cdots  &   f_{2n}(x)\\
    0   &   0   &   f(x)    &   \cdots  &   f_{3n}(x)\\
    \vdots  &   \vdots  &   \vdots  &   \ddots  &   \vdots\\
    0   &   0   &   0   &   \cdots  &   f(x)\\
\end{bmatrix}\]

We have three cases:
\begin{enumerate}
    \item Suppose that $g(x) \neq 0 \neq f(x)$. Then looking at the diagonal entries of $F(x)G(x)F(x) = 0$ and we see that $f(x)g(x)f(x) = 0$. Now, these are polynomials in $R[x]\setminus\{0\}$. We use the fact that $R$ is inner McCoy to replace $g(x)$ with a nonzero $d\in R$. Let $D := dE_{1n}$. Then $F(x)DF(x) = 0$. 

\item Suppose $f(x) \neq 0 = g(x)$. Now, we know that some entries of $G(x)$  are nonzero due to the choice of $G(x)$. Suppose that the $k^{th}$ row of the matrix $G(x)$ is the last row with nonzero entries. If the leading nonzero entry of this row appears at the $m^{th}$ column, looking at the $(k,m)^{th}$ entry of the product $F(x)G(x)F(x) = 0$ gives us the expression $f(x)g_{km}(x)f(x) = 0$, where $g_{km}(x)$ is the nonzero $(k,m)^{th}$ entry of $G(x)$.
As $R$ is an inner McCoy ring, we can replace $g_{km}(x)$ with a nonzero $d\in R$ such that $f(x)df(x) = 0$.
Let $D:=dE_{km}$. Then $F(x)DF(x) = 0$.

\item Suppose that $f(x) = 0$.Then it is easy to see that $ F(x)E_{1n}F(x) = (F(x)E_{11})E_{1n}F(x)  = 0$. 
\end{enumerate}

\noindent
This exhausts all possible cases and so $DUT_n(R)$ is inner McCoy. This holds for all $n\geq 2$.
\\
Now, suppose that $DUT_n(R)$ is inner McCoy, for $n\geq 2$ and let $g(x), f(x) \in R[x]\setminus\{0\}$ such that $f(x)g(x)f(x) = 0$. If $G(x):= g(x)I_n$, $F(x):=f(x)I_n$, then we have $F(x)G(x)F(x) = 0$. As $DUT_n(R)$ is inner McCoy, we have a nonzero $D \in DUT_n(R)$ such that $F(x)DF(x) = 0$. If 

\[ D = 
\begin{bmatrix}
    d & d_{12} &  d_{13}  & \cdots & d_{1n}\\
    0 & d & d_{23} & \cdots & d_{2n} \\
    0 & 0 & d & \cdots & d_{3n}\\
    \vdots & \vdots & \vdots & \ddots & \vdots \\
    0 & 0 & 0 & \cdots & d
\end{bmatrix}
\]

then we have two cases here:
\begin{enumerate}
    \item If $d \neq 0$, then the diagonal entries of $F(x)DF(x) = 0$ gives us $f(x)df(x) = 0$, which is what we  require.
    \item If $d =0$, then as $D \neq 0$, there are some entries that are nonzero. Suppose the $k^{th}$ row of $D$ is the last row with nonzero entries. Then if the leading nonzero entry of the $k^{th}$ row of $D$ is in the $m^{th}$ column, we have that the $(k,m)^{th}$ entry of the matrix $F(x)DF(x) = 0$ gives $f(x)d_{km}f(x) = 0$, which gives us a nonzero inner annihilator of $f(x)$ in $R$.
\end{enumerate}
This exhausts all cases and so $R$ must be inner McCoy. 
\end{proof}

If we were able to construct a McCoy/inner McCoy ring with only trivial idempotents, we would be able to use this result to construct more rings that are both \emph{i}-reversible and McCoy/inner McCoy. In fact, we shall discuss such an example below.

\begin{eg}\label{ex2}
    Consider the ring, 
    $L = \mathbb{Z}_2 \langle c_0, c_1, d_0, d_1 \rangle$        and the ideal,
        \begin{equation*}
            J = \langle c_0d_0, c_0d_1 + c_1d_0, c_1d_1, d_id_j (0 \leq i, j \leq), d_ic_j (0 \leq i,j \leq 1) \rangle.
        \end{equation*}
        The quotient ring $S := L/J$ is not McCoy, as Nielsen showed that $S$ is left McCoy but not right McCoy, as the polynomials $F(x) = \Bar{c_0} +\Bar{c_1}x$ and $G(x) = \Bar{d_0}+\Bar{d_1}x \in S[x]$ are such that $F(x)G(x)=0$ but $F(x)$ can never be annihilated on the right by a nonzero element of $R$ \cite{nielsen2006}. %It is also clear that the ring is not semi-commutative, as 
        %\begin{equation*}
         %   c_0d_0 = 0 \neq c_0c_1d_0.
       % \end{equation*}
        However, we claim that this ring is indeed inner McCoy.
%\end{eg}
%\begin{proof}
        First, we make the observation that a general element in $R$ is of the form
    \begin{equation}\label{genelt}
    \begin{split}
    \gamma & = f_0 + f_1(c_0) + f_2(c_0)d_1 +f_3(c_0)f_4(c_1) + f_5(c_0)f_6(c_1)d_0
    \end{split}
    \end{equation}
    where $f_0 \in \mathbb{Z}_2$ and $f_1(x), f_2(x), f_3(x), f_4(x), f_5(x), f_6(x) \in \mathbb{Z}_2[x]$. This is possible by the use of Bergman's diamond lemma \cite{bergman1978} to obtain a basis for the given ring in terms of the indeterminates. Note that in $\gamma$, each of the indeterminates $a_i$ and $b_i$ are not elements of the ring $K$, but of the quotient ring $R$. We are abusing notation here for brevity. 

    Let $P(x)$ and $Q(x)$ be nonzero polynomials in $R[x]$ such that $Q(x)P(x)Q(x) = 0.$ Suppose $P(x) = \displaystyle\sum_{i=0}^{m}p_ix^i$ and $Q(x) = \displaystyle\sum_{i=0}^{n}q_ix^i$.
    If each $q_i$ has a zero constant term, then $d_0Q(x) = 0 \Rightarrow Q(x)d_0Q(x) = 0.$
    We now show that for any polynomial $Q(x) \in R[x]$ whose coefficients have nonzero constant terms can never be inner annihilated by a nonzero polynomial.
    
    Let $k \in \{0, 1, 2, \cdots, n\}$ be the smallest index such that $q_k \neq 0$. 
    We now define for each coefficient $p_i$ of $P(x)$ a corresponding element $p'_i \in R$ as follows:
    \[ p'_i = \begin{cases}
    0 & \textrm{if } p_i =0\\
    \sum \alpha_i & \textrm{if } p_i \neq 0
    \end{cases}
    \]
     where $\alpha_i$'s are all the terms of smallest degree in $p_i \in R$. Then we select $j \in \{0,1,2, \cdots, m\}$ such that $p'_j$ is the smallest amongst all the $p_i$'s defined above. We assert the existence of such $j$ as $P(x) = 0$. Now, in the degree $j+2k$ term of the expression $Q(x)P(x)Q(x) = 0$, we have
    \[\sum_{\substack{r,s,t: \\ r+s+t = j+2k}} q_rp_sq_t = 0 \in R.\]   
    
    As $I$ is a homogeneous ideal, we have that all like degree terms of any fixed degree above must add to zero. But from choice of $j$ and $k$, we have the term $p'_j \cdot 1$, which is nonzero and cannot cancel with any other term in the summation above. We have arrived at a contradiction here and thus, each $q_i$ must have a zero constant term and we are done.
\end{eg}

\begin{lem} \label{c3}
    The quotient ring $R$ as specified in Example \ref{ex2}, is a ring with only trivial idempotents.
\end{lem}
\begin{proof}
    We will suppose to the contrary that this ring does have a nontrivial idempotent, say $e$. The $e$ will be as given in equation (\ref{genelt}). We make a few observations regarding this ring:
    \begin{enumerate}
        \item The idempotent $e$ has a corresponding orthogonal idempotent $1-e$. As the ring has characteristic 2 this is just $1+e$.
        \item Looking at each monomial of this ring, we notice that, with the exception of the monomials $1$, $c_0^{i}$ and $c_0^{i}c_1^{j}$ where $i, j \in \mathbb{N}$, all the monomials are nilpotent with nilpotency 2.
    \end{enumerate}
    We will, without loss of generality, make the assumption, that $f_0 = 0$ and so the constant term in $1+e$ is $1$. As $e$ is nontrivial, we have some monomials of least nonzero degree in $e$. We will group all such non-zero monomials of least degree in $e$, which we will call $e_0$. We will group all the other terms together also and label them $\alpha$.\\
    So we have,
    \begin{equation*}
        (1+e)e =        (1+e_0 +\alpha)(e_0 + \alpha) =
        e_0 + \alpha +e_0^2 +e_0\alpha +\alpha e_0 + \alpha^2  = 0.
    \end{equation*}
    We make the observation that all monomials here all evaluate to zero or have degree strictly greater than those of $e_0$. Now, the ideal $I$ is homogeneous, and so the terms of same degree in the above expression must add to zero. But the term $e_0$ is nonzero, making $(1+e)e \neq 0$, which gives us a contradiction to $e$ being an idempotent. Our assumption was therefore incorrect, and so the ring has only trivial idempotents.
\end{proof}

In view of Theorem \ref{c1}, Lemma \ref{c3} and Proposition \ref{c2}, we can provide an example of a ring that is not reversible, is \emph{i}-reversible, is not McCoy and is inner McCoy.

\begin{eg}\label{ex3}
    The rings $DUT_n(R)$, where $R$ is the quotient ring in Example \ref{ex2} and $n \geq 2$, are rings that are not reversible, are \textit{i}-reversible, are not McCoy and are inner McCoy.
\end{eg}

\begin{remark}
Similar to Example \ref{ex3}, we can choose $R$ appropriately to construct rings which are inner McCoy and not \textit{i}-reversible and rings which are neither \textit{i}-reversible nor inner McCoy. The ring $M_2(\mathbb R)$ is not inner McCoy as seen in Example \ref{mn} and it is \emph{i}-reversible, from \cite[Theorem 4.3]{khurana2020}. 
\end{remark}

\begin{remark}
    The second half of the proof of Theorem \ref{c1} may be replicated to show that if $T_n(R)$ is inner McCoy then $R$ is inner McCoy. The converse however is not true. For example if we let $R$ be the ring from Example \ref{ex2} and consider $F(x)=\begin{bmatrix}
        c_1x+c_0&1\\
        0&x+1
    \end{bmatrix}$ and $G(x)=\begin{bmatrix}
        d_1x+d_0&0\\ 0&0
    \end{bmatrix}$ then it is clear that $F(x)G(x)F(x)=0$ but only the zero matrix will inner annihilate $F(x)$.
\end{remark}

While the classes of i-reversible and inner McCoy rings are primarily independent of each other we may ask about relationships between inner McCoy rings and semi-commutative rings. The triangular matrix ring $T_2(R)$ is inner McCoy but not semi-commutative for $E_{11}E_{22}=0$ but $E_{11}\begin{bmatrix}
    1&1\\ 0&1
\end{bmatrix}E_{22}=E_{12}\neq0$. Also, the ring $DUT_2(\mathbb R)$ is inner McCoy (see Theorem \ref{c1}) as well as semi-commutative for if $AB=0$ in $DUT_2(\mathbb R)$ then the diagonal entries of $A$ and $B$ must be 0 and it is easy to see that $A\cdot DUT_2(\mathbb R)\cdot B=0$.

Clearly, there are inner McCoy rings that are semi-commutative and those that are not semi-commutative. It still remains to be seen if there are examples of rings that are semi-commutative but not inner McCoy. One could therefore ask the following:

\textbf{Question:} Is every semi-commutative ring inner McCoy? 

We can summarise all the relationships between all the known classes of rings and the class of inner McCoy rings as follows:
\begin{enumerate}
    \item All McCoy rings are inner McCoy. In particular, Camillo and Nielsen \cite{camillo2008} have already mentioned that reduced and  Armendariz rings are McCoy, so they are inner McCoy as well.
    \item There are left/right McCoy rings that are inner McCoy but not McCoy.
    \item There are rings that are inner McCoy and semi-commutative as well as rings that are inner McCoy but not semi-commutative.
    \item There are rings that are inner McCoy, and are \textit{i}-reversible, that are inner McCoy and are not \emph{i}-reversible, that are not inner McCoy but are \emph{i}-reversible, and that are neither inner McCoy nor \emph{i}-reversible.
\end{enumerate}

\section{Maximal inner McCoy subring of the full matrix ring}
It is already well known in literature \cite{camillo2008} that even if $R$ is a McCoy ring, then neither the full matrix ring $M_n(R)$ nor the upper triangular matrix ring $T_n(R)$ need be McCoy. We have already seen that over a commutative ring $R$, the full matrix ring $M_n(R)$ is not inner McCoy while  the triangular matrix rings $T_n(R)$ and $DUT_n(R)$ are inner McCoy. We now ask if $M_n(R)$ has a maximal inner McCoy subring, that is a subring that does have the inner McCoy property such that no subring larger than it is inner McCoy. First, we note that for inner McCoy rings, the maximal inner McCoy subring will be the ring itself. In what follows we will work with matrix rings over a field $F$.

\begin{theorem}
    For all $n \geq 2$, the ring $T_n(F)$ is the maximal inner McCoy subring of $M_n(F)$.
\end{theorem}
\begin{proof}
If $T_n(F)$ is not a maximal inner McCoy subring of $M_n(F)$, then we should be able to find subrings of $M_n(F)$ that properly contain $T_n(F)$. The only such rings are the generalised upper triangular matrices over full matrix rings $M_{r_1}(\mathbb{R}), M_{r_2}(\mathbb{R}), \cdots M_{r_k}(\mathbb{R})$, where $\displaystyle\sum_{i=1}^{k} r_i = n$ and $r_i > 1$ for all $1 \leq i \leq k$. It is enough to look at rings of the form: $$R_r = \begin{bmatrix}
    M_r(F) & M_{r \times (n-r)}(F) \\
    \mathbb O_{(n-r)\times r} & M_{n-r}(F)
\end{bmatrix}$$   
where $r = 1,2,\cdots, \lfloor \frac{n}{2} \rfloor $ upto conjugation. Suppose $$A_r(x): = \begin{bmatrix}
    x & x^2 &  \cdots & x^r \\
    x^{r+1} & x^{r+2} & \cdots & x^{2r} \\
    \vdots & \vdots & \ddots & \vdots \\
    x^{r^2-r+1} & x^{r^2-r+2} & \cdots & x^{r^2}
\end{bmatrix}, B_r(x) := \begin{bmatrix}
    -x^{r-1} & -x^{r-1} & \cdots & -x^{r-1} \\
        x^{r-2} & x^{r-2} & \cdots & x^{r-2} \\
        \vdots & \vdots & \ddots & \vdots \\
        1 & 1 & \cdots & 1
\end{bmatrix}.$$ We know that $A_r(x)B_r(x)A_r(x) = 0$. We construct the following elements in $R_r[x]$:
$$A(x):= \begin{bmatrix}
    A_r(x) & \mathbb O_{r\times (n-r)} \\
    \mathbb O_{(n-r)\times r} & x^{r^2}A_{n-r}(x)
\end{bmatrix}, B(x):= \begin{bmatrix}
    B_r(x) & \mathbb O_{r\times (n-r)} \\
    \mathbb O_{(n-r)\times r} & \mathbb O_{(n-r)\times (n-r)}
\end{bmatrix}.$$
It is clear that in $R_r[x]$, $A(x)B(x)A(x) = 0$. If we try to replace $B(x)$ by a matrix $M \in R_r$, it is straightforward to see that all the entries of $M$ must be zero in order to ensure $A(x)MA(x) = 0$. This shows that none of the $R_r$'s are inner McCoy. One can show this is the same for any generalised upper triangular matrix ring as specified before in a similar way. Therefore $T_n(F)$ is a maximal inner McCoy subring of $M_n(F)$.
\end{proof}

Theorem \ref{c1} gives us the relationship between inner McCoyness of $R$ and the ring extension $DUT_n(R)$ while Theorem \ref{x1} gives us a condition for when the $T_n(R)$ is inner McCoy. We now ask about other ring extensions. We discuss this in the next section.

\section{Some Extensions of Inner McCoy Rings}

In this section we discuss some results on creating new inner McCoy rings from older ones. In particular, we look at polynomial rings (in one and several variables) and direct products of rings.

\begin{theorem}
    A ring $R$ is inner McCoy if and only if the polynomial ring $R[x]$ is inner McCoy.
\end{theorem}
\begin{proof}
    We assume that the ring $R$ is inner McCoy. 
    Let $F(y) = \displaystyle\sum_{i=0}^{n}f_i y^i$ and $G(y) = \displaystyle\sum_{j=0}^{m} g_j y^j$
    be two nonzero polynomials from $(R[x])[y]\setminus \{0\}$ such that $G(y)F(y)G(y) = 0$.
    \par
    We know that the coefficients of both $F(y)$ and $G(y)$ are of the form,
    $f_i = \displaystyle\sum_{s=0}^{p_i}a_{is}x^s$
    and $g_j = \displaystyle\sum_{t=0}^{q_j}b_{jt}x^t$.
    We define the number $k$ as follows,
    \begin{equation*}
        k :=\sum_{i=0}^n deg_x(f_i) + 2\sum_{j=0}^m deg_x(g_j)
    \end{equation*}
    where we have $deg_x(0):=0$. Now, we can consider the statement $G(y)F(y)G(y) =0$.
    This expression holds for any formal commuting symbol in place of $y$. In particular, if we replace $y$ with $x^k$, we get $G(x^k)F(x^k)G(x^k) =0$.
    We choose this particular power of $x$ so that all the coefficients of this statement match with that in the original expression that we were dealing with.
    \par
    We may think of this more formally as the image of a homomorphism from $(R[x])[y]$ to $R[x]$, that maps $y$ to $x^k$, and preserves coefficients for this particular value of $k$.
    \par
    However, this expression is one in $R[x]$. As we know that $R$ is inner McCoy, we have a nonzero $u \in R$ such that
    \begin{equation*}
        G(x^k)uG(x^k) = 0.
    \end{equation*}
    Again, as $x^k$ is a formal symbol, we return it as $x^k = y$ and obtain an inner annihilator of $G(y)$ in $R \subset R[x]$. This proves that $R[x]$ is inner McCoy.\\
    Now, we suppose that $R[x]$ is inner McCoy.\par
    Let 
    $
        f(x) = \displaystyle\sum_{i=0}^{n} a_ix^i
    $
    and
    $
        g(x) = \displaystyle\sum_{j=0}^{m} b_jx^j
    $
    be a pair of nonzero polynomials of $R[x]$ such that
    $
        g(x)f(x)g(x) = 0.
    $ \par
    We can simply change the formal symbol to being $y$ and this does not affect the validity of the above expression. So, we have
    $
        g(y)f(y)g(y) =0.
    $
    However, this is an expression in $(R[x])[y]$. Therefore, we can use the assumption that $R[x]$ is inner McCoy to obtain a polynomial
    \begin{equation*}
        h(x) = \sum_{k=0}^{p} h_kx^k \in R[x]\setminus \{0\}
    \end{equation*}
    such that
    $
        g(y)h(x)g(y) = 0.
    $
    In particular, we have 
    $
        g(y)h_pg(y) =0 
    $
    and we are assured that $h_p$, being the leading coefficient of $h(x)$ is nonzero. This means that $g(y)$ is inner annihilated by a nonzero element of $R$.
    \\
    Now, as the formal symbol $y$ does not play any special role as is, the above expression still holds if we replace $y$ with $x$. So, we have
    $
        g(x)h_pg(x) =0 
    $
    in $R[x]$. This shows that $R$ is inner McCoy.
\end{proof}
This parallels a similar result provided by Lei, Chen and Ying \cite{lei2007} that does the same for McCoy rings. We now obtain a corollary to this theorem that seems natural:

\begin{cor}
    Suppose that $I$ is an indexing set and $\{x_\alpha \}_{\alpha \in I}$ are a set of commuting indeterminates. Then we can say that the ring $R$ is inner McCoy if and only if the ring $R[\{x_\alpha\}_{\alpha \in I}]$ is inner McCoy.
\end{cor}
\begin{proof}
    We first note that if we have $I$ to be any finite set, then the result follows by an induction argument. \\
    Suppose that $I$ is any nonfinite set. Then we make the observation that for any pair of nonzero polynomials $f, g \in R[\{x_\alpha \}_{\alpha \in I}]$, we can only have a finite number of indeterminates in each polynomial.     We take the subring $S$ of $R[\{x_\alpha \}_{\alpha \in I}]$ of as many finite indeterminates such that both $f$ and $g \in S$, and use the above argument to conclude the same. We can repeat this for any pair of nonzero polynomials in $R[\{x_\alpha \}_{\alpha \in I}]$ where one inner annihilates the other to get the result required.
\end{proof}

\begin{theorem}
    Let $S$ be a multiplicatively closed subset consisting of only central regular elements of an inner McCoy ring $R$, then $S^{-1}R$ is inner McCoy.
\end{theorem}
\begin{proof}
    Let $f(x)=\displaystyle\sum_{i=0}^n\frac{a_i}{t_i}x^i$, $g(x)=\displaystyle\sum_{j=0}^m\frac{b_j}{s_j}x^j\in S^{-1}R$ such that $f(x)g(x)f(x)=0$. Let $w_i=\displaystyle\prod_{\substack{k=0\\k\neq i}}^nt_k\ $,$v_j=\displaystyle\prod_{\substack{k=0\\k\neq j}}^ms_k$ and let $f_1(x)=\displaystyle\sum_{i=0}^na_iw_ix^i\ $, $g_1(x)=\displaystyle\sum_{j=0}^mb_jv_jx^j$. Then $f_1(x)g_1(x)f_1(x)=0$ in $R$. By the inner McCoy condition we can find a non-zero $c\in R$ such that $f_1(x)cf_1(c)=0$. Then it is easy to see that $\dfrac{c}{1}$ will inner annihilate $f(x)$.
\end{proof}

\begin{cor}
    The Laurent polynomial ring $R[x,x^{-1}]$ of an inner McCoy ring is inner McCoy.
\end{cor}

We now conclude with a very simple, yet natural question about how we can construct inner McCoy rings from known ones. We present here another such way involving the direct product of rings;

\begin{theorem}\label{520}
    Suppose $R_1$ and $R_2$ are inner McCoy rings. Then their direct product $R_1 \times R_2$ is also inner McCoy.
\end{theorem}

\begin{proof}
    Consider the nonzero polynomials of $(R_1\times R_2)[x]$ as follows:
    $
        f(x) = \displaystyle\sum_{i=0}^{m} (a_i, b_i)x^i
    $
    and 
    $
        g(x) = \displaystyle\sum_{j=0}^{n} (c_j, d_j)x^j
    $
    where 
    $
        g(x)f(x)g(x) = 0
    $.
    We will then consider the polynomials that correspond to each summand of $R_1 \times R_2$ as:
    \begin{equation*}
        f_1(x) = \sum_{i=0}^m a_ix^i,
    \quad 
        f_2(x) = \sum_{i=0}^m b_ix^i,
    \quad
        g_1(x) = \sum_{j=0}^n c_jx^j,
    \quad
        g_2(x) = \sum_{j=0}^m d_jx^j.
    \end{equation*}
    With these we can state that for each summand, we obtain:
    \begin{equation*}
        g_1(x)f_1(x)g_1(x) = 0
    \quad \textrm{and} \quad 
        g_2(x)f_2(x)g_2(x) = 0.
    \end{equation*}
    However, as both summands $R_1$ and $R_2$ are each inner McCoy, we can find a nonzero element in each ring $r_1 \in R_1$ and $r_2 \in R_2$ respectively such that
    \begin{equation*}
        g_1(x)r_1g_1(x) = 0 \quad
    \textrm{and} \quad
        g_2(x)r_2g_2(x) = 0.
    \end{equation*}
    Then $(r_1, r_2) \in R_1\times R_2 \setminus \{0\}$ is the required nonzero inner annihilator of $g(x)$ from $ (R_1\times R_2)$
\end{proof}

We make note of the fact that this theorem parallels a similar result by \cite{ying2007} for McCoy rings. We extend this result by induction to a countable number of rings that are all inner McCoy.

\begin{cor}
    Let  $R_1, R_2, \cdots R_n$ be a collection of $n$ inner McCoy rings. Then the direct product of these rings is also inner McCoy.
\end{cor}
\begin{proof}
    This is easy to show by using an induction argument and using Theorem \ref{520} as a base case. 
\end{proof}

With this result, we can actually give an example of a ring that is neither left McCoy nor right McCoy, but is inner McCoy.

\begin{eg}
    Let $R$ be the ring as specified in Example \ref{ex2}. We know that $R$ is left McCoy but not right McCoy, and we have shown that $R$ is inner McCoy. Also, taking the opposite ring $R^{opp}$ gives us an example of a ring that is right McCoy, not left McCoy and is inner McCoy.
    The ring $S:= R \times R^{opp}$  is  neither left nor right McCoy, but is inner McCoy.
    
    To show $S$ is not right McCoy, we take the polynomial $f(x)=(\Bar{a_0}, 0)+(\Bar{a_1},0)x+(\Bar{a_2},0)x^2+(\Bar{a_3},0)x^3$ and $g(x) = (\Bar{b_0},0)+(\Bar{b_1},0)x$, such that $f(x)g(x) = 0$, but there is no nonzero $r = (r_1, r_2) \in S$ such that $f(x)r = 0$. 
    
    To show $S$ is not left McCoy, we take the polynomial $f(x) = (0, \Bar{a_0})+(0, \Bar{a_1})x+(0, \Bar{a_2})x^2+(0, \Bar{a_3})x^3$ and $g(x) = (0, \Bar{b_0})+(0, \Bar{b_1})x$ such that $g(x)f(x) = 0$, but there is no nonzero $r = (r_1, r_2) \in S$ such that $rf(x) = 0$.

    The ring $S$ is  neither left nor right McCoy, but as both summands are inner McCoy, $S$ must be inner McCoy using Theorem \ref{520}.
\end{eg}

Finally, we present an example of a ring that is not McCoy but is both semi-commutative and inner McCoy.

\begin{eg}\label{ex1}
    Let $K := \mathbb{Z}_2\langle a_0, a_1, a_2, a_3, b_0, b_1 \rangle$, the free associative algebra with unity over $\mathbb{Z}_2$. 
    Let $I$ be the ideal generated as follows:
    \begin{multline*}
        I   = \langle a_0b_0, a_0b_1+a_1b_0, a_1b_1+a_2b_0, a_2b_1 +a_3b_0, a_3b_1, \\
    a_0a_j (0 \leq j \leq 3), a_3a_j (0 \leq j \leq 3), a_1a_j+a_2a_j (0 \leq j \leq 3), \\
    b_ib_j (0 \leq i,j \leq 1), b_ia_j (0 \leq i \leq 1,0 \leq j \leq 3) \rangle.
    \end{multline*}
    The ring $R := \dfrac{K}{I}$ is semi-commutative and left McCoy, but not right McCoy, as the polynomials $F(x) = \Bar{a_0} +\Bar{a_1}x+\Bar{a_2}x^2 +\Bar{a_3}x^3$ and $G(x) = \Bar{b_0}+\Bar{b_1}x \in R[x]$ are such that $F(x)G(x)=0$ but $F(x)$ can never be annihilated on the right by a nonzero element of $R$ \cite{nielsen2006}. Therefore, this ring is not McCoy. However, we note that this ring is actually inner McCoy. This can be shown by repeating the arguments presented in example \ref{ex2}.
%\end{eg}
\end{eg}

\bibliography{refs.bib}

\begin{thebibliography}{10}

\bibitem{mccoy1942}
N.~H. McCoy, ``Remarks on divisors of zero,'' {\em The American Mathematical
  Monthly}, vol.~49, no.~5, pp.~286--295, 1942.

\bibitem{hong2000}
C.~Y. Hong, N.~K. Kim, and T.~K. Kwak, ``Ore extensions of baer and pp-rings,''
  {\em Journal of Pure and Applied Algebra}, vol.~151, no.~3, pp.~215--226,
  2000.

\bibitem{hong2010}
C.~Y. Hong, N.~K. Kim, and Y.~Lee, ``Extensions of {M}c{C}oy's theorem,'' {\em
  Glasgow Mathematical Journal}, vol.~52, no.~1, pp.~155--159, 2010.

\bibitem{weiner1952}
L.~G. Jones and L.~Weiner, ``Problem 4419,'' {\em The American Mathematical
  Monthly}, vol.~59, no.~5, pp.~336--337, 1952.

\bibitem{rege1997armendariz}
M.~B. Rege and S.~Chhawchharia, ``Armendariz rings,'' {\em Proceedings of the
  {J}apan academy, series {A}, mathematical sciences}, vol.~73, pp.~14--17,
  1997.

\bibitem{nielsen2006}
P.~P. Nielsen, ``Semi-commutativity and the {M}c{C}oy condition,'' {\em Journal
  of Algebra}, vol.~298, no.~1, pp.~134--141, 2006.

\bibitem{lam2014}
T.~Lam and P.~P. Nielsen, ``Inner inverses and inner annihilators in rings,''
  {\em Journal of Algebra}, vol.~397, pp.~91--110, 2014.

\bibitem{jung2019}
D.~W. Jung, C.~I. Lee, Y.~Lee, S.~Park, S.~J. Ryu, H.~J. Sung, and S.~J. Yun,
  ``On reversibility related to idempotents,'' {\em Bulletin of the Korean
  Mathematical Society}, vol.~56, no.~4, pp.~993--1006, 2019.

\bibitem{zhao2009}
R.~Zhao and Z.~Liu, ``Extensions of {M}c{C}oy rings,'' in {\em Algebra
  Colloquium}, vol.~16, pp.~495--502, World Scientific, 2009.

\bibitem{camillo2008}
V.~Camillo and P.~P. Nielsen, ``Mc{C}oy rings and zero-divisors,'' {\em Journal
  of Pure and Applied Algebra}, vol.~212, no.~3, pp.~599--615, 2008.

\bibitem{lei2007}
Z.~Lei, J.~Chen, and Z.~Ying, ``A question on {M}c{C}oy rings,'' {\em Bulletin
  of the Australian Mathematical Society}, vol.~76, no.~1, pp.~137--141, 2007.

\bibitem{lama2022}
V.~B. Lama, B.~N. Suhas, S.~Mazumdar, and R.~DSouza, ``Extensions of
  \emph{i}-reversible rings,'' {\em Journal of Algebra and its Applications},
  2022.

\bibitem{bergman1978}
G.~M. Bergman, ``The diamond lemma for ring theory,'' {\em Advances in
  mathematics}, vol.~29, no.~2, pp.~178--218, 1978.

\bibitem{khurana2020}
A.~Khurana and D.~Khurana, ``\emph{i}-{R}eversible rings,'' {\em Journal of
  Algebra and its Applications}, vol.~19, no.~04, p.~2050076, 2020.

\bibitem{ying2007}
Z.~Ying, J.~Chen, and Z.~Lei, ``Extensions of {M}c{C}oy rings,'' {\em
  Northeastern Mathematical Journal}, vol.~24, no.~1, pp.~85--94, 2008.

\end{thebibliography}
\bibliographystyle{ieeetr}

\par

\par

\par

\end{document}